\documentclass[a4paper,11pt]{article}

\input{styleart.sty}

\title{Elementary subgroups of the free group are free factors - a new proof}
\date{\today}
\author{Chlo\'{e} Perin}

\begin{document}

\maketitle

\begin{abstract} 
In this note we give a new proof of the fact that an elementary subgroup (in the sense of first-order theory) of a non abelian free group $\F$ must be a free factor. The proof is based on definability of orbits of elements of under automorphisms of $\F$ fixing a large enough subset of $\F$.  
\end{abstract}

\section{Introduction} 

In 1945, Tarski asked whether non abelian free groups of different (finite) ranks are elementary equivalent, that is, whether they satisfy exactly the same first-order sentences in the language of groups. At the turn of the millennium, the question was finally given a positive answer (see \cite{Sel6}, as well as \cite{KharlampovichMyasnikov}). In fact, the proofs showed something stronger, namely that for any $2\leq k\leq n$, the canonical embedding of the free group $\F_k$ of rank $k$ in the free group $\F_n$ of rank $n$ is an elementary embedding. The converse was proved in \cite{PerinElementary}:
\begin{theorem} \label{ElementaryIsFreeFactor} Let $H$ be a subgroup of a nonabelian free group $\F$. If the embedding of $H$ in $\F$ is elementary, then $H$ is a free factor of $\F$.
\end{theorem}

In a work in progress, Guirardel and Levitt give another proof of this result, using homogeneity of the free group (proved in \cite{OuldHoucineHomogeneity} and \cite{PerinSklinosHomogeneity} independently).  

In this note, we give yet another proof of Theorem \ref{ElementaryIsFreeFactor}, which relies on the (highly non trivial) fact proved in \cite{PerinSklinosForking} that if $A$ is a large enough subset of parameters in a free group $\F$ then orbits under the group of automorphisms of $\F$ fixing $A$ pointwise are definable by a first-order formula with parameters in $A$. This relative definability of orbits holds also for concrete torsion free hyperbolic groups - namely freely indecomposable torsion free hyperbolic groups which do not admit a structure of extended hyperbolic tower over a proper subgroup (for example, groups whose cyclic JSJ decomposition does not have any surface type vertices are concrete). Hence we also get a proof that concrete torsion free hyperbolic groups do not admit proper elementary subgroups (see Theorem \ref{ConcreteNoElementary}).  

Recall that a subgroup $H$ of a group $G$ is elementary in the sense of first-order logic if the tuples of elements of $H$ have the same type over $H$ and over $G$, that is, if they have the same first-order properties when seen as tuples in $H$ and as tuples in $G$. More precisely, $H$ is elementary in $G$ if
given any formula $\phi(x_1, \ldots, x_n)$ in the language of groups with free variables $x_1, \ldots, x_n$ and any elements $h_1, \ldots, h_n$ of $H$, the formula $\phi(h_1, \ldots, h_n)$ holds in $H$ if and only if it holds in $G$.

Note that it is straightforward to prove
\begin{lemma} \label{ElementaryIsFG} Let $\F$ be a finitely generated free group. Let $H$ be an elementary subgroup of $\F$. Then $H$ is finitely generated. In fact, it must be a retract of $\F$.
\end{lemma}
(A subgroup $G'$ of a group $G$ is called a retract of $G$ if there exists a surjective homomorphism $G \to G'$ which restricts to the identity on $G'$).

However, free groups admit retracts which are not free factors (see \cite{TurnerTestWords} - for example, the subgroup generated by $a[a,b]$ in the free group on $\{a,b\}$ is a retract but the element $a[a,b]$ is not primitive). Therefore, the lemma is not enough to prove Theorem \ref{ElementaryIsFreeFactor} and the tools used in the proof of definability of orbits over large sets of parameters - on which we rely to prove Theorem \ref{ElementaryIsFreeFactor} - are much more advanced.

Still, let us give a proof of Lemma \ref{ElementaryIsFG} as a warm up.
\begin{proof} Let $a_1, \ldots, a_n$ be a basis for $\F$. Since $H$ is also free, we can choose $\{h_1, h_2, \ldots \}$ to be a (possibly infinite) basis for $H$.
Each $h_i$ can be written as a word $w_i(a_1, \ldots, a_n)$, thus for any $k$ which is at most the rank of $H$, the sentence 
$$ \psi_k(h_1, \ldots, h_k): \exists x_1, \ldots, x_n \bigwedge^{k}_{i=1} h_i = w_i(x_1, \ldots, x_n)$$
holds in $\F$.

Suppose that the rank of $H$ is at least $n+1$. The formula $\psi_{n+1}(h_1, \ldots, h_{n+1})$ holds in $\F$ so by assumption, it must hold in $H$. In particular there exist elements $b_1, \ldots, b_n$ of $H$ which generate a subgroup $B$ of $H$ containing the free factor $H_1 = \langle h_1, \ldots, h_{n+1} \rangle$ of $H$. By Grushko's theorem, $B \cap H_1 = H_1$ is a free factor of $B$, which is a contradiction since $H_1$ has rank $n+1$ while $B$ is generated by $n$ elements. Thus $H$ has rank at most $n$.

Now consider the sentence $\psi_n(h_1, \ldots, h_n)$: it holds in $\F$ hence by assumption it holds in $H$, that is, there exist elements $b_1, \ldots, b_n$ of $H$ such that $h_i=w_i(b_1, \ldots, b_n)$ for all $i \in \{ 1, \ldots, n \}$. The morphism $\theta: \F \to H$ defined by $\theta(a_j)= b_j$ satisfies $\theta(h_i) = \theta(w_i(a_1, \ldots, a_n)) = w_i(b_1, \ldots, b_n)=h_i$, so it is surjective and it restricts to the identity on $H$. 
\end{proof}

\section{Relative JSJ decompositions and automorphism groups}

Let $G$ be a torsion free hyperbolic group which is freely indecomposable relative to some non trivial subgroup $H$. Denote by $\Aut_H(G)$ the group of automorphisms of $G$ fixing $H$ pointwise.

\paragraph{Modular group.} The modular group of $G$ relative to $H$, which we denote $\Mod_H(G)$, is the subgroup of $\Aut_H(G)$ generated by Dehn twists associated to one edge splittings of $G$ over cyclic groups relative to $H$, i.e. splitting of $G$ as an amalgamated product $G = A*_Z B$ or as an HNN extension $G = A*_Z$ for which $Z$ is infinite cyclic and $H \leq A$.

It is a result of Rips and Sela (see \cite{RipsSelaHypI}) that (in the case of a freely indecomposable torsion free hyperbolic group), the modular group $\Mod_H(G)$ has finite index in the group of automorphisms $\Aut_H(G)$.

\paragraph{JSJ decomposition.} We assume familiarity with Bass-Serre theory (see \cite{SerreTrees}). A JSJ-decomposition for a group $G$ relative to a subgroup $H$ over a class ${\cal A}$ of subgroups is a graph of groups decomposition $\Lambda$ for $G$ which is universal in the sense that it "encodes" (with some precise meaning we do not define here) all possible splittings of $G$ as a graph of groups whose edge groups belong to the class ${\cal A}$ and in which $H$ is elliptic.  

It is a consequence of \cite{RipsSelaJSJ}, \cite{Bowditch} that in the case where $G$ is torsion free hyperbolic and freely indecomposable relative to $H$, and ${\cal A}$ is the class of cyclic subgroups of $G$, such a decomposition $\Lambda$ exists. It has infinite cyclic edge groups, and $H$ is elliptic in $\Lambda$ - in other words we can assume it is contained in a vertex group $G_H$ of $\Lambda$.

The universal property of JSJ decompositions implies in particular that for any splitting of $G$ as an amalgamated product $G = A*_Z B$ or as an HNN extension $G = A*_Z$ for which $Z$ is infinite cyclic and $H \leq A$, the vertex group $G_H$ of $\Lambda$ containing $H$ is contained in $A$.

Moreover, it is possible to choose $\Lambda$ to be characteristic in the following sense: the action of $G$ on the tree $T$ corresponding to $\Lambda$ is invariant under precomposition by automorphisms of $G$ fixing $H$ pointwise, that is, for any element $\theta \in \Aut_H(G)$ there exists an automorphism $j: T \to T$ such that for any $x \in T$ and $g \in G$ we have $j(g \cdot x) = \theta(g) \cdot j(x)$. The tree associated with this decomposition is sometimes called the tree of cylinders of the JSJ deformation space (see \cite{GuirardelLevittUltimateJSJ}). We call this decomposition $\Lambda$ the \textbf{canonical cyclic JSJ decomposition of $G$ relative to $H$}. 

For a more detailed summary and pointers to the relevant literature see Section 4 of \cite{PerinSklinosForking}.

\begin{lemma} \label{AutFixBaseVertexGp} Let $G$ be a torsion free hyperbolic group which is freely indecomposable relative to some non abelian subgroup $H$. Let $\Lambda$ be the canonical cyclic JSJ decomposition for $G$ relative to $H$, and let $G_H$ denote the vertex group containing $H$. 

Then the orbit of any element of $G_H$ under $\Aut_H(G)$ is contained in $G_H$.
\end{lemma}

\begin{proof} Denote by $T$ the tree associated to $\Lambda$. Note that there is a unique vertex $v$ stabilized by $H$ since $H$ is not cyclic. Let $G_H$ denote its stabilizer.

Let $\theta \in \Aut_H(G)$, and consider the  automorphism $j: T \to T$ such that for any $x \in T$ and $g \in G$ we have $j(g \cdot x) = \theta(g) \cdot j(x)$. 

The group $\theta(G_H)$ is the stabilizer of the vertex $j(v)$, but $H \leq \theta(G_H)$ hence $H$ stabilizes both $v$ and $j(v)$ - we must have $v=j(v)$, and $\theta(G_H) = G_H$.
\end{proof}

We also prove
\begin{lemma} \label{OrbitsOfBaseGpFinite} Let $G$ be a torsion free hyperbolic group which is freely indecomposable relative to some non abelian subgroup $H$. Let $\Lambda$ be the canonical cyclic JSJ decomposition for $G$ relative to $H$, and let $G_H$ denote the vertex group containing $H$. 

Then the orbit of any element of $G_H$ under $\Aut_H(G)$ is finite.
\end{lemma}

\begin{proof} We show in fact that elements of $G_H$ are fixed under $\Mod_H(G)$. Since $\Mod_H(G)$ has finite index in $\Aut_H(G)$ the result follows easily. 

In any splitting of $G$ as an amalgamated product $G = A*_Z B$ or as an HNN extension $G = A*_Z$ for which $Z$ is infinite cyclic and $H \leq A$, the vertex group $G_H$ of $\Lambda$ containing $H$ is contained in $A$. This implies that Dehn twists induced by such splittings restrict to a conjugation on $G_H$. Now such a Dehn twist restricts to the identity on $H$ if and only if it is in fact the identity on $A$ (centralizers of non abelian subgroups of torsion free hyperbolic groups are trivial). Thus it must be in fact the identity on $G_H$ as well. Since $\Mod_H(G)$ is generated by such Dehn twists, this remains true of any modular automorphism. 
\end{proof}

\section{Definability of orbits over big sets of parameters} 

The following is Theorem 5.3 of \cite{PerinSklinosForking}:
\begin{theorem}  \label{DefinabilityOrbitsConcrete} If $G$ is a torsion-free hyperbolic group and $A$ a subset of $G$ such that $G$ is concrete relative to $A$, then orbits of tuples of $G$ under $\Aut_A(G)$ are definable over $A$. 
\end{theorem}
To say $G$ is concrete over $A$ is to say that $G$ is freely indecomposable with respect to $A$ and admits no proper structure of extended hyperbolic tower over $A$ - see Definition 5.1 of \cite{PerinSklinosForking}. 

The following is Lemma 5.19 of \cite{PerinElementary}:
\begin{lemma} \label{FreeGroupsNoHypFloor} If $\F$ is a free group which is freely indecomposable relative to a subset $A$ then $\F$ does not admit a proper structure of hyperbolic floor with respect to $A$.
\end{lemma}
In particular, $\F$ is in that case concrete relative to $A$. Thus we have
\begin{corollary} \label{DefinabilityOrbits} Let $\F$ be a non abelian free group. Let $A \subseteq \F$ be such that $\F$ is freely indecomposable with respect to $A$. Then the orbit of any tuple $g$ of $\F$ under automorphisms of $\F$ fixing $A$ pointwise is definable over $A$.
\end{corollary}

In fact, in the case of the free group we will need a slightly more general statement, namely
\begin{theorem} \label{DefOrbitsCor} Let $\F$ be a non abelian free group. Let $A \subseteq \F$, and let $\F_A$ be the smallest free factor of $\F$ containing $A$. Then for any tuple $g$ contained in $\F_A$, there exists a formula with parameters in $A$, say $\psi_g(x, A)$, such that $\F \models \psi_g(x, A)$ if and only if there exists an automorphism $\theta$ of $\F$ fixing $A$ pointwise such that $\theta(g) =x$.
\end{theorem}

The proof is very similar to that of Theorem 5.3 of \cite{PerinSklinosForking}.
\begin{proof} By Lemma 3.7 in \cite{PerinSklinosHomogeneity}, there is a finite subset $A_0$ of $A$ such that $\F_A$ is freely indecomposable with respect to $A_0$, $A$ is elliptic in any $\F_A$-tree in which $A_0$ is elliptic, and 
$\Mod_{A_0}(\F_A) = \Mod_{A}(\F_A)$. By Corollary 4.5 in \cite{PerinSklinosHomogeneity}, we can assume moreover 
that any embedding $j: \F_A \to \F_A$ which restricts to the identity on $A_0$ restricts to the identity on $A$.

By Theorem 4.4 in \cite{PerinSklinosHomogeneity}, there is a finite set of quotients $\{\eta_j: \F_A \to Q_j\}_{j=1, \ldots, m}$ 
such that any non injective endomorphism $\theta: \F_A \to \F$ which restricts to the identity on $A_0$ factors through one 
of the quotient maps $\eta_j$ after precomposition by an element $\sigma$ of $\Mod_A(\F_A)$. 
For each $j$, choose $u_j$ a non trivial element in $\Ker(\eta_j)$. 

Let $\gamma_1, \ldots , \gamma_r$ be a generating set for $\F_A$. Write each element $a$ of $A_0$  as a
word $w_a(\gamma_1, \ldots, \gamma_r)$, each element $u_j$ a word $w_{u_j}(\gamma_1, \ldots, \gamma_r)$, and the tuple $g$ as a tuple of words $w_{g}(\gamma_1, \ldots, \gamma_r)$.

Let $\Lambda$ be a JSJ decomposition of $\F_A$ with respect to $A$. Two morphisms $h$ and $h'$ from $\F_A$ to $\F$ are said to be 
$\Lambda$-related if $h$ and $h'$ coincide up to conjugation on rigid vertex groups of $\Lambda$, and for any flexible 
vertex group $S$ of $\Lambda$, $h(S)$ is non abelian if and only if $h'(S)$ is non abelian. It is easy to see that there is a 
formula  $\Rel(\bar{x}, \bar{y})$ such that for any pair of morphisms $h$ and $h'$ from $\F_A$ to $\F$, the morphism 
$h'$ is $\Lambda$-related to $h$ if and only if $\F \models \Rel(h(\gamma_1, \ldots, \gamma_r), h'(\gamma_1, \ldots, \gamma_r))$ 
(see Lemma 5.18 in \cite{ThesisPerin}).

Consider now the following formula $\phi(z, A_0)$:
\begin{eqnarray*}
 \exists x_1, \ldots, x_r \; && \left\{ z = w_g(x_1, \ldots, x_r) \wedge \bigwedge_{a \in A_0} a=w_a(x_1, \ldots, x_r) \right\} \\
&& \wedge \; \forall y_1, \ldots, y_r \left \{ \Rel(\bar{x},\bar{y}) \rightarrow \bigvee_j w_{u_j}(y_1, \ldots, y_r) \neq 1 \right\}.
\end{eqnarray*}

Suppose $\F \models \phi(g', A_0)$. Then the morphism $h: \F_A \to \F$ given by $\gamma_j \mapsto x_j$ sends 
$g$ to $g'$ and fixes $A_0$, moreover no endomorphism $h'$ which is $\Lambda$-related to $h$ factors through one of the maps $\eta_i$. 
This implies that $h$ is injective. By our choice of $A_0$, we get that $h$ fixes $A$ pointwise. Moreover, $h(\F_A)$ is freely indecomposable relative to $h(A)=A$ - this means we must have $h(\F_A) \leq \F_A$.  But by the relative co-Hopf property for torsion-free hyperbolic groups 
(see Corollary 4.2 of \cite{PerinSklinosHomogeneity}), this in turn implies that $h(\F_A) =\F_A$. Thus it is easy to extend $h$ to an automorphism of $\F$ fixing $A$. Thus the set defined by the formula $\phi(z, A_0)$ is contained in the orbit of $g$ by $\Aut_A(\F)$.

To finish the proof it is enough to show that $\F\models \phi(g,A_0)$. 

It is obvious that the first part 
of the sentence is satisfied by $\F$ (just take $x_j$ to be $\gamma_j$). If the second part is not satisfied with this choice of $x_j$, this 
means that there exists a morphism $h': \F_A \to \F$ which is $\Lambda$-related to the embedding of $\F_A$ in $\F$, and which kills 
one of the elements $u_j$. Thus $h'$ restricts to conjugation on the rigid vertex groups of $\Lambda$, 
sends surface type flexible vertex groups on non abelian images, and is non injective: 
it is a non injective preretraction $\F_A \to \F$ with respect to $\Lambda$. By Proposition 5.11 in \cite{PerinElementary}, there exists a non injective preretraction $\F_A \to \F_A$ with respect to $\Lambda$. By Proposition 5.11 in \cite{PerinElementary}, 
this implies that $\F_A$ admits a structure of hyperbolic floor over $A$, a contradiction to Lemma \ref{FreeGroupsNoHypFloor}. 
\end{proof}

\section{A new proof of Theorem \ref{ElementaryIsFreeFactor}} 

We can now give the proof
\begin{proof}
Let $\F$ be a non abelian free group, and let $H$ be an elementary subgroup of $\F$. In particular, $H$ is not abelian. Denote by $\F_H$ the smallest free factor of $\F$ containing $H$, consider the JSJ decomposition $\Lambda$ of $\F_H$ relative to $H$, and denote by $G_H$ the vertex group containing $H$. 

Suppose first that $\Lambda$ is not trivial - then there exists $g \in \F_H$ which does not lie in $G_H$. Note that by Lemma \ref{AutFixBaseVertexGp}, the orbit of $g$ under $Aut_H(\F_H)$ does not meet $G_H$, so in particular it does not meet $H$. By Theorem \ref{DefOrbitsCor}, there is a formula $\psi_g(x, H)$ defining the orbit of $g$ under $\Aut_H(\F)$ in $\F$, that is, $\F \models \psi_g(x, H)$ iff there exists $\theta \in \Aut_H(\F)$ such that $\theta(g)=x$.

Now the formula $\exists x \psi_g(x, H)$ holds in $\F$ so it should hold in $H$. This gives an element $g_0 \in H$ so that $H \models \psi_g(g_0, H)$, but then we must have $\F \models \psi_g(g_0, H)$ contradicting the fact that the orbit of $g$ does not meet $H$.

If $\Lambda$ is trivial, but $H \neq \F_H$, take an element $g \in \F_H -H$. Its orbit under $\Aut_H(\F)$ is finite by Corollary \ref{OrbitsOfBaseGpFinite}, say equal to $\{g=g_1, \ldots, g_k\}$, and definable in $\F$ by the formula $\psi_g(x,H)$. We can write the following formula
$$ \exists x_1, \ldots , x_k \bigwedge^k_{i=1} \psi_g(x_i,H) \wedge \bigwedge_{i \neq j} x_i\neq x_j $$
which holds in $\F$ and says that the set defined in $\F$ by $\psi_g(x, H)$ has at least $k$ elements. 

Now the elements satisfying $\psi_g(x,H)$ in $H$ satisfy $\psi_g(x,H)$ in $\F$ as well, hence they are exactly the elements of $H$ which lie in the orbit of $g$ under $\Aut_H(\F)$. But in $H$ there are strictly less than $k$ elements of the orbit $\{g=g_1, \ldots, g_k\}$ of $g$ since $g \not \in H$. The formula above does not hold in $H$: this is a contradiction.
\end{proof}

In fact, using Theorem \ref{DefinabilityOrbitsConcrete}, we can apply similar ideas to torsion free hyperbolic groups which are concrete relative to any proper subgroup, i.e. which are freely indecomposable and do not admit any structure of hyperbolic floor. Such groups are simply called \textbf{concrete}. We get
\begin{theorem} \label{ConcreteNoElementary}
Let $G$ be a concrete torsion free hyperbolic group. Then $G$ has no proper elementary subgroups. 
\end{theorem}
 
Indeed, if $G$ is torsion free hyperbolic and abelian it is trivial or infinite cyclic and the result is well known. If $G$ is not abelian, and $H$ is an elementary subgroup of $G$ then on the one hand $H$ is non abelian and on the other hand, by definition, $G$ is concrete relative to $H$. Thus by Theorem \ref{DefinabilityOrbitsConcrete}, the orbits of tuples of elements of $G$ under $\Aut_H(G)$ are definable over $H$. By following the proof of Theorem \ref{ElementaryIsFreeFactor}, we get the result.

\bibliography{biblio}
\end{document}